\documentclass[11pt]{article}
\usepackage[table]{xcolor}
\usepackage{graphicx}
\usepackage{amsmath,amsfonts,amssymb}
\usepackage{graphicx}
\usepackage{multirow}
\usepackage{tikz,pgf}
\usepackage{url}
\usepackage{amsmath}
\usepackage{graphicx}
\usepackage{epsfig}
\usepackage{epstopdf}
\usepackage{color}
\usepackage{enumitem}
\def\tto{\;{\lower 1pt \hbox{$\rightarrow$}}\kern -10pt
\hbox{\raise 2pt \hbox{$\rightarrow$}}\;}

\def\Bar{\overline}
\def\ra{\rangle}
\def\la{\langle}
\def\ve{\varepsilon}
\def\epsilon{\varepsilon}

\def\h{\hfill\Box}
\def\R{\Bbb R}

\def\N{\Bbb N}
\def\ox{\bar{x}}

\def\oy{\bar{y}}

\def\ri{\mbox{\rm ri}}
\def\icr{\mbox{\rm icr}}

\def\gph{\mbox{\rm gph}}
\def\aff{\mbox{\rm aff}}
\def\epi{\mbox{\rm epi}}
\def\hypo{\mbox{\rm hypo}}

\def\dom{\mbox{\rm dom}}
\def\aff{\mbox{\rm aff}}

\def\sint{\mbox{\rm int}}

\def\cone{\mbox{\rm cone}}

\def\iri{\mbox{\rm iri}\,}

\def\h{\hfill\square}

\def\emp{\emptyset}
\def\st{\stackrel}

\def\lm{\lambda}
\def\gg{\gamma}

\def\al{\alpha}

\def\emp{\emptyset}
\def\st{\stackrel}

\def\lm{\lambda}
\def\gg{\gamma}

\def\al{\alpha}

\def\qri{\mbox{\rm qri}}
\setlist[enumerate,1]{itemsep=0.0ex,parsep=0.5ex,label={\rm(\alph*)},leftmargin=*, align=left}
\newcounter{lk}

\usepackage[colorlinks=true]{hyperref}

\begin{document}
\begin{center}
{\sc\bf Fenchel-Rockafellar Theorem in Infinite Dimensions \\via Generalized Relative Interiors}\\[1ex]
{\sc D. V. Cuong} \footnote{Department of Mathematics, Faculty of Natural Sciences, Duy Tan University, Da Nang, Vietnam (dvcuong@duytan.edu.vn). This research is funded by Vietnam National Foundation for Science and Technology Development (NAFOSTED) under grant number 101.02-2020.20}$^,$\footnote{American Degree Program, Duy Tan University, Da Nang, Vietnam.}, {\sc B. S. Mordukhovich}\footnote{Department of Mathematics, Wayne State University, Detroit, Michigan 48202, USA (boris@math.wayne.edu). Research of this author was partly supported by the USA National Science Foundation under grants DMS-1512846 and DMS-1808978, by the USA Air Force Office of Scientific Research grant \#15RT04, and by Australian Research Council under grant DP-190100555.},
{\sc N. M.  Nam}\footnote{Fariborz Maseeh Department of Mathematics and Statistics, Portland State University, Portland, OR
97207, USA (mnn3@pdx.edu).},
{\sc G. Sandine}\footnote{Fariborz Maseeh Department of
 Mathematics and Statistics, Portland State University, Portland, OR 97207, USA (gsandine@pdx.edu  ).}\\[2ex]
 {\bf Dedicated to Phan Quoc Khanh on the occasion of his 75th birthday}
\end{center}
\small{\bf Abstract.} In this paper we provide further studies of Fenchel duality theory in the general framework of locally convex topological vector (LCTV) spaces. We prove the validity of  the  Fenchel strong  duality under some qualification conditions via generalized relative interiors  imposed on the epigraphs and the domains of the functions involved. Our results directly generalize the classical Fenchel-Rockafellar theorem on strong duality from finite dimensions to LCTV spaces. \\[1ex]

\newtheorem{Theorem}{Theorem}[section]
\newtheorem{Proposition}[Theorem]{Proposition}
\newtheorem{Remark}[Theorem]{Remark}
\newtheorem{Lemma}[Theorem]{Lemma}
\newtheorem{Corollary}[Theorem]{Corollary}
\newtheorem{Definition}[Theorem]{Definition}
\newtheorem{Example}[Theorem]{Example}
\renewcommand{\theequation}{\thesection.\arabic{equation}}
\normalsize
\section{Introduction}
\setcounter{equation}{0}

Duality theory has a central role in optimization theory and its applications. From a \emph{primal optimization problem} with the objective function defined in a primal space, a  \emph{dual problem}  in the dual space is formulated with the hope that the new problem is easier to solve, while having a close relationship with the primal problem.   Its important role in optimization has made duality theory attractive for extensive research over the past few decades; see, e.g., \cite{Anh-Khanh2016,B-S-S,Bauschke2011,bl,Borwein2000,bot1,bot-revis,bot2,bot3,brezis,daniele2,durea,AG,HU1,HU2,ktz,mnrt,NgS,zduality} and the references therein.

Given a proper convex function $f$ and a proper concave function $g$ defined on $\R^n$, the classical Fenchel-Rockafellar theorem asserts that \begin{equation}\label{Eq-Fen-Roc} \inf\{f(x)-g(x)\
|\ x\in \R^n\}= \sup\{g_*(x^*)-f^*(x^*)\ |\ x^*\in \R^n\}, \end{equation} under the assumption that $\ri(\dom(f))\cap \ri(\dom(g))\neq \emptyset$; see  \cite[Theorem~31.1]{r}. In the setting of
an LCTV space, the validity of \eqref{Eq-Fen-Roc} requires the continuity assumption of either $f$ or $g$ at some point in the intersection of their domains; see, e.g.,
\cite[Theorem~1.12]{brezis}. Nevertheless, this assumption does not reduce to the classical relative interior condition in  \cite[Theorem~31.1]{r}. To overcome this shortcoming, generalized
interior concepts were introduced and applied to further study this theorem. Among many generalized interior concepts, we refer the readers to the notion of \emph{quasi-relative interior}
introduced by Borwein and Lewis in \cite{bl}. Under the assumption that $A(\qri(\dom g))\cap \ri(\dom h) \neq\emptyset$, Borwein and Lewis proved the following strong duality theorem:
\begin{equation*}\label{Eq-Fen-Roc1}
  \inf\{g(x)+h(Ax)\ |\ x\in X\}= \sup\{-g^*(A^*y^*)-h^*(-y^*)\ |\ y^*\in Y^*\},
\end{equation*} where $g\colon X\to (-\infty,\infty]$ defined on an LCTV space $X$ and $h\colon \R^n\to (-\infty,\infty]$ are proper convex functions, and $A\colon X\to \R^n$ is a continuous
linear mapping with its adjoint $A^*$; see \cite[Corollary~4.3]{bl}. Note that this new result reduces to the classical Fenchel-Rockafellar theorem in the finite-dimensional case with the use of
the identity mapping $A$, but it does not provide a full generalization to infinite dimensions.  Further significant success in this research direction was achieved by  Bo\c{t}, Z\u{a}linescu,
and others; see \cite{bot1,bot2,bot3,NgS,zduality, zonuseqri} and the references therein.

The  main goal of this paper is to provide a full generalization of the Fenchel-Rockafellar theorem to the LCTV setting using a number of known  generalized relative  interior concepts  along
with the regularity condition introduced in our recent paper \cite{CBN}. As a consequence, we obtain a
fully infinite-dimensional generalization of the Fenchel-Rockafellar theorem for functions defined on LCTV spaces.  Our paper is organized as follows. In Section 2 we introduce basic definitions and  clarify several important known
results involving generalized relative interiors used throughout the paper. Section~3 is devoted to the study of generalized relative interiors for convex graphs of set-valued mappings between
LCTV spaces. In this section, we will extend the results of \cite[Corollary 9(iii)]{zonuseqri} and \cite[Theorem 5.6]{CBN} to set-valued mappings and provide a similar result for mappings between
LCTV spaces when the codomain is partially ordered. Section 4 is devoted to the study of the Fenchel-Rockafellar theorem in LCTV spaces.  Throughout the paper, we use standard notation and results
from convex analysis which can be found in \cite{z}. All spaces under consideration are Hausdorff real topological vector spaces unless otherwise stated.

\section{On Generalized Relative Interiors} \setcounter{equation}{0} In this section we introduce some standard notation and definitions from convex analysis in topological vector spaces used
throughout the paper. The reader can find more details in the book by Z\u{a}linescu \cite{z}.  We also provide a survey with detailed clarifications of some known facts involving the generalized
relative interiors of convex sets.   Throughout this paper, consider a real Hausdorff topological vector space $X$ with its topological dual $X^*$ and the canonical pairing $\la x^*,x\ra:=x^*(x)$
with $x\in X$ and $x^*\in X^*$.  For a nonempty subset $\Omega$ of $X$, $\Omega$ is a {\em cone} if $tw\in\Omega$ whenever $w\in\Omega$ and $t\ge0$, and the
{\em conic hull} of $\Omega$ is defined by $\cone(\Omega):=\big\{tw\in X\;|\;t\ge 0,\;w\in \Omega\big\}$.

Given a convex subset $\Omega$ of a topological vector space $X$,
the {\em relative interior} of $\Omega$ is defined by
\begin{equation}\label{ri} \mbox{\rm
ri}(\Omega):=\big\{x\in\Omega\;\big|\;\exists\;\mbox{\rm
neighborhood}\,V\,\mbox{\rm of }\,x\,\mbox{\rm with
}\,V\cap\Bar{\aff}(\Omega)\subset\Omega\big\}. \end{equation} If $X$
is finite-dimensional, this notion reduces to the classical one
formulated, e.g., in \cite[Definition 1.68]{bmn}, since  the
closure operation is not needed in \eqref{ri} due to the automatic
closedness of affine sets in finite dimensions.\vspace*{0.03in}

We begin with deriving some useful characterizations of relative interiors of convex sets in $\R^n$  that are important for the subsequent extensions of this notion to convex sets in infinite dimensions.
For the reader's convenience, we present here the definition of the {\em normal cone} to a convex subset $\Omega$ of a topological vector space at a point $\ox\in\Omega$:
\begin{equation}\label{ncd0}
N(\ox;\Omega):=\big\{x^*\in X^*\;\big|\;\la x^*,x-\ox\ra\le 0\;\mbox{ for all }\;x\in \Omega\big\}\;\mbox{ if }\;\ox\in\Omega.
\end{equation}

\begin{Theorem}\label{pts2} Let $\Omega$ be a nonempty convex set in $\R^n$, and let $\ox\in\R^n$. The following properties are equivalent:
\begin{enumerate}
\item $\ox\in\ri(\Omega)$.
\item $\ox\in\Omega$ and for every $x\in\Omega$ with $x\ne\ox$ there exists a vector $u\in\Omega$ such that $\ox\in(x,u)$.
\item $\ox\in\Omega$ and $\cone(\Omega-\ox)$ is a linear subspace of $\R^n$.
\item $\ox\in\Omega$ and $\overline{\cone}(\Omega-\ox)$ is a linear subspace of $\R^n$.
\item $\ox\in\Omega$ and the normal cone $N(\ox;\Omega)$ is a subspace of $\R^n$.
\end{enumerate}
\end{Theorem}
The obtained finite-dimensional characterizations of relative interior motivate the major extensions of this notion to infinite dimensions, which are considered in what follows.

\begin{Definition}\label{qri} Let $\Omega$ be a convex subset of a topological vector space $X$. Then we have:\\[1ex]
{\rm(a)} The {\sc intrinsic relative interior} of $\Omega$ is the set
\begin{equation*}\label{iri}
\mbox{\rm iri}(\Omega):=\big\{x\in\Omega\;\big|\;\mbox{\rm cone}(\Omega-x)\;\mbox{\rm is a subspace of }\;X\big\}.
\end{equation*}
{\rm(b)} The {\sc quasi-relative interior} of $\Omega$ is the set
\begin{equation*}\label{qri1}
\mbox{\rm qri}(\Omega):=\big\{x\in\Omega\;\big|\;\Bar{\mbox{\rm cone}}(\Omega-x)\;\mbox{\rm is a subspace of }\;X\big\}.
\end{equation*}
{\rm(c)} We say that a convex set $\Omega\subset X$ is {\sc quasi-regular} if $\qri(\Omega)={\rm iri}(\Omega)$.
\end{Definition}

Due to Theorem~\ref{pts2}, both notions in Definition~\ref{qri}(a,b)
reduce to the relative interior of $\Omega$ in finite-dimensional
spaces. The one in (a) is also known under the name ``intrinsic
core" of $\Omega$, which may be confusing. Definition~\ref{qri}(c)
designates the property $\qri(\Omega)={\rm iri}(\Omega)$ by labeling
the sets satisfying this condition as {\em quasi-regular} ones. The
latter property plays an important role in the subsequent results of
this section. Some sufficient conditions for the quasi-regularity
property of convex sets are presented below.\vspace*{0.03in}

To proceed further, we first present a simple equivalent description of intrinsic relative interior points of nonempty convex sets (see: \cite[Lemma 2.3]{BG}). Recall that a point $\ox\in\Omega$ is
{\em relatively absorbing} for $\Omega$ if for every $x\in\Omega\setminus\{\ox\}$ there exists $u\in\Omega$ such that $\ox\in(x,u)$.

\begin{Proposition}\label{iri-desc} Let $\Omega$ be a nonempty convex subset of a topological vector space $X$, and let $\ox\in\Omega$. Then we have $\ox\in\iri(\Omega)$ if and only if $\ox$ is a
relatively absorbing point of the set $\Omega$.
\end{Proposition}
\noindent {\bf Proof.}  To verify the ``if" part, observe that relatively absorbing points of $\Omega$ can be equivalently described as follows: for any $x\in\Omega\setminus\{\ox\}$ there exists $\al>1$
such that $(1-\al) x+\al\ox\in\Omega$. Pick now any $v\in\cone(\Omega-\ox)$ (with $v\neq 0$) and find $\lm>0$ with $\lm v+\ox\in\Omega$. The relative absorbability of $\ox$ gives us $\al>1$ such that
\begin{equation*}
(1-\al)(\lm v+\ox)+\al\ox=\lm(\al-1)(-v)+\ox\in\Omega.
\end{equation*}
This yields $-v\in\cone(\Omega-\ox)$, and hence $\cone(\Omega-\ox)$ is a subspace of $X$.

To justify the converse implication, pick any $\ox\in\iri(\Omega)$ and $x\in\Omega$ (with $x\neq \ox$). Since $\cone(\Omega-\ox)$ is a subspace of $X$,
\begin{equation*}
x-\ox\in\cone(\Omega-\ox)\;\mbox{ and }\;\ox-x\in\cone(\Omega-\ox).
\end{equation*}
Choose $t>0$ such that $\ox-x=t(w-\ox)$ with $w\in \Omega$. Thus we get some number $\al>1$ ($\al=1+1/t$) for which
\begin{equation*}
(\al-1)(\ox-x)+\ox=(1-\al)x+\al\ox\in\Omega.
\end{equation*}
The latter means that $\ox$ is a relatively absorbing point of $\Omega$. $\h$

Note that Proposition~\ref{iri-desc} shows that the equivalence (b)$\Longleftrightarrow$(c) of Theorem~\ref{pts2} holds for nonempty convex sets in $X$. Furthermore, the other equivalence
(d)$\Longleftrightarrow$(e) of Theorem~\ref{pts2} suggests a similar relationship in the general LCTV setting established in what follows.

\begin{Definition} Let $X$ be a topological vector space.
\begin{enumerate}
\item Let $\Omega$ be a subset of  $X$. Define the {\sc polar} of $\Omega$ by
\begin{equation*}
\Omega^\circ:=\{x^*\in X^*\; |\; \la x^*, w\ra\leq 1\; \mbox{\rm for all }w\in \Omega\}.
\end{equation*}
\item Let $\Theta$ be a subset of $X^*$. Define the {\sc polar} of $\Theta$ by
\begin{equation*}
\Theta^\circ:=\{x\in X\; |\; \la z^*, x\ra\leq 1\; \mbox{\rm for all }z^*\in \Theta\}.
\end{equation*}
\end{enumerate}
\end{Definition}
It follows from the definition that if $\Omega$ is a cone in $X$, then
\begin{equation*}
\Omega^\circ:=\{x^*\in X^*\; |\; \la x^*, w\ra\leq 0\; \mbox{\rm for all }w\in \Omega\}.
\end{equation*}
Similarly, if $\Theta$ is a cone in  $X^*$, then
\begin{equation*}
\Theta^\circ:=\{x\in X\; |\; \la z^*, x\ra\leq 0\; \mbox{\rm for all }z^*\in \Theta\}.
\end{equation*}

\begin{Lemma} \label{polar relationship} Let $C$ be a nonempty convex cone in an  LCTV space $X$. Then \begin{equation*} (C^\circ)^\circ=\Bar{C}. \end{equation*} \end{Lemma} \noindent {\bf Proof.}
Observe that $C^\circ$ is a nonempty convex cone in $X^*$. It follows directly from the definition that $(C^\circ)^\circ$ is a closed subset of $X$ and $C\subset (C^\circ)^\circ$. Then
$\Bar{C}\subset (C^\circ)^\circ$. Now, fix any $x\in (C^\circ)^\circ$ and suppose on the contrary that $x\notin \Bar{C}$. By the convex strict separation theorem, there exists $x^*\in X^*$ such
that \begin{equation*} \la x^*, u\ra \leq 0\; \mbox{\rm for all }u\in C\; \mbox{\rm and }\la x^*, x\ra>0. \end{equation*} Thus $x^*\in C^\circ$ and $\la x^*, x\ra>0$ which contradicts the fact
that $x\in (C^\circ)^\circ$ and completes the proof. $\h$

\begin{Lemma}\label{normal cone and polar} Let $\Omega$ be a nonempty convex set in a topological vector space $X$, and let $\ox\in \Omega$. Then
\begin{equation*}
N(\ox; \Omega)=\Theta^\circ=\Bar{\Theta}^\circ,
\end{equation*}
where $\Theta:=\mbox{\rm cone}(\Omega-\ox)$.
\end{Lemma}
\noindent {\bf Proof.} Fix any $x^*\in N(\ox; \Omega)$. It follows from the definition that
\begin{equation*}
\la x^*, x-\ox\ra\leq 0\; \mbox{\rm for all }x\in \Omega.
\end{equation*}
This yields $\la x^*, w\ra\leq 0$ for all $w\in \Theta$, and hence $x^*\in \Theta^\circ$. For the reverse inclusion, taking any $x^*\in \Theta^\circ$ gives that $\la x^*, w\ra\leq 0$ for all $w\in \Theta$. Then for any
$x\in \Omega$ we have  $x-\ox\in \Theta$, and so  $\la x^*, x-\ox\ra\leq 0$. This implies $x^*\in N(\ox; \Omega)$, and thus $N(\ox; \Omega)=\Theta^\circ$. It is also straightforward to show that
$\Theta^\circ=\Bar{\Theta}^\circ$. $\h$

The following result known by \cite[Proposition 2.8]{bl} presents a simple equivalent description of quasi-relative interior points of nonempty convex sets.
\begin{Proposition}\label{qrinormal}  Let $\Omega$ be a nonempty convex subset of an LCTV space $X$, and let $\ox\in\Omega$. Then we have
\begin{equation}\label{qri-char}
\big[\ox\in\qri(\Omega)\big]\Longleftrightarrow\big[N(\ox;\Omega)\;\mbox{ is a linear subspace of }\;X^*\big]
\end{equation}
\end{Proposition}
\noindent {\bf Proof.} Suppose first that $\ox\in\qri(\Omega)$. It follows from the definition that the set $\Bar{\Theta}$, where $\Theta:=\mbox{\rm cone}(\Omega-\ox)$, is a linear subspace of $X$. An easy exercise shows that
$\Bar{\Theta}^\circ$ is also a linear subspace of $X^*$. Then Lemma \ref{normal cone and polar} tells us that $N(\ox; \Omega)=\Theta^\circ$ is a linear subspace of $X^*$. Conversely, suppose $N(\ox; \Omega)$ is a linear subspace of $X^*$. Using Lemma \ref{polar relationship} and Lemma \ref{normal cone and polar}, we have
\begin{equation*}
N(\ox; \Omega)^\circ=(\Bar{\Theta}^\circ)^\circ=\Bar{\Theta}.
\end{equation*}
Since $N(\ox; \Omega)$ is a linear subspace of $X^*$, the set $N(\ox; \Omega)^\circ$ is also a linear subspace of $X$. Thus, $\Bar{\Theta}$ is a linear subspace of $X$, and therefore $\ox\in \qri(\Omega)$. $\h$

We continue with a well-known version of {\em proper separation} of a singleton from a convex set that gives us yet another characterization of quasi-relative interior; see: \cite[Theorem 2.3]{Flores1}.

\begin{Proposition}\label{qri-sep1} Let $\Omega$ be a convex set in an LCTV space $X$, and let $\ox\in\Omega$. Then the sets $\{\ox\}$ and $\Omega$ are properly separated if and only if
$\ox\notin\qri(\Omega)$.
\end{Proposition}
\noindent {\bf Proof.} Using the normal cone characterization \eqref{qri-char} of quasi-relative interior points, we get that $\ox\notin\mbox{\rm qri}(\Omega)$ if and only if there exists
$x^*\in N(\ox;\Omega)$ with $-x^*\notin N(\ox;\Omega)$. It follows from the normal cone construction \eqref{ncd0} for convex sets that $\la x^*,x\ra\le\la x^*,\ox\ra$ for all $x\in\Omega$. Then the
inclusion $-x^*\notin N(\ox;\Omega)$ gives us an $x_0\in\Omega$ such that $\la-x^*,x_0\ra>\la-x^*,\ox\ra$, which is equivalent to $\la x^*,x_0\ra<\la x^*,\ox\ra$ and hence justifies the statement of the
proposition. $\h$

Next we establish relationships between the notions of relative, intrinsic relative, and quasi-relative interiors of  convex sets in LCTV spaces and give conditions that ensure quasi-regularity of a set.
These results can be found in \cite{BG,bl,zduality}.
\begin{Theorem}\label{ri-rel} Let $\Omega$ be a  convex subset of a topological vector space $X$. Then we have the inclusions
\begin{equation}\label{ri-rel1}
\ri(\Omega)\subset\iri(\Omega)\subset\qri(\Omega).
\end{equation}
If furthermore $X$ is locally convex and  $\ri(\Omega)\ne\emp$, then the inclusions in \eqref{ri-rel1} become the equalities
\begin{equation}\label{eq-ri-quari}
\ri(\Omega)=\iri(\Omega)=\qri(\Omega).
\end{equation}
\end{Theorem}
\noindent {\bf Proof.} We first show that $\ri(\Omega)\subset\iri(\Omega)$. Take $\ox\in\ri(\Omega)$ and fix $x\in\Omega$ with $x\ne\ox$. It follows from \eqref{ri} that $\ox\in\Omega$ and there exists a
neighborhood $V$ of $\bar x$ such that
\begin{equation}\label{affc0}
V\cap\overline{\aff}(\Omega)\subset\Omega.
\end{equation}
Choose $0<t<1$ so small that $u:=\ox+t(\ox-x)\in V$. Then $u\in\aff(\Omega)\subset\overline{\aff}(\Omega)$, and we get from \eqref{affc0} that $u\in\Omega$. It follows that
\begin{equation*}
\ox=\frac{t}{1+t}x+\frac{1}{1+t}u\in(x,u),
\end{equation*}
which therefore verifies by Proposition~\ref{iri-desc} that $\bar{x}\in\iri(\Omega)$. This tells us that $\ri(\Omega)\subset\iri(\Omega)$. The other inclusion $\iri(\Omega)\subset\qri(\Omega)$ in
\eqref{ri-rel1} is trivial, since the subspace property of $\cone(\Omega-\bar{x})$ clearly implies that the closure $\overline{\cone}(\Omega-\bar{x})$ is also a linear subspace of $X$.

To prove the equalities in \eqref{eq-ri-quari}, it is sufficient to show that if $\ri(\Omega)\ne\emp$ and $\bar{x}\in\qri(\Omega)$, then $\bar{x}\in\ri(\Omega)$. Arguing by contradiction, assume that
$\bar{x}\notin\ri(\Omega)$ and begin with the case where $\ox=0$. If $0\notin\overline{\Omega}$, then by the strict separation there exists $x^*\in X^*$ such that
\begin{equation}\label{eq-ri-qri2}
\la x^*,x\ra<0\;\text{ for all }\;x\in\Omega.
\end{equation}
In the complement setting where  $0\in\overline{\Omega}\setminus\ri(\Omega)$, denote $X_0:=\overline{\aff}(\Omega)$ and get $0\in X_0$ telling us that $X_0$ is a closed subspace of $X$. It is easy to see
that $0\notin\ri(\Omega)=\sint_{X_0}(\Omega)$, where $\sint_{X_0}(\Omega)$ is the interior of $\Omega$ with respect to the space $X_0$. Applying the separation result to the sets $\Omega$ and $\{0\}$ in
the topological  space $X_0$, we find $x_0^*\in X_0^*$ ensuring that
\begin{equation}\label{eq-stric-sepa}
\la x^*_0,x\ra\leq 0\;\text{ for all }\;x\in\Omega,\; \mbox{\rm and }  \la x^*_0, \bar w\ra<0\; \mbox{\rm for some }\bar w\in \Omega.
\end{equation}

Then the Hahn-Banach extension theorem from \cite[Theorem 2.10]{Bonnans2000} shows that there exists an extension $x^*\in X^*$ of $x_0^*$ such that
\begin{equation*}
\la x^*,x\ra\le 0\;\mbox{ for all }\;x\in\Omega.
\end{equation*}
In either case there exists $x^*\in X^*$ such that $\la x^*, x\ra \leq 0$ for all $x\in \Omega$ and hence for all $x\in\overline{\cone}(\Omega)$.  Since $0\in\qri(\Omega)$, we have that
$\overline{\cone}(\Omega)$ is a linear subspace, and therefore
\begin{equation*}
\la x^*,x\ra=0\;\mbox{ for all }\;x\in\overline{\cone}(\Omega).
\end{equation*}
This contradicts the conditions in \eqref{eq-stric-sepa} and also in \eqref{eq-ri-qri2}, and thus verifies that $0\in\ri(\Omega)$. Turning finally to the general case for $\ox$, we reduce it to the case
where $\ox=0$ due to the obvious relationships
\begin{equation*}
\begin{aligned}
&\bar{x}\in\ri(\Omega)\Longleftrightarrow 0\in\ri(\Omega-\bar{x})=\ri(\Omega)-\{\bar{x}\}\;\mbox{ and}\\
&\bar{x}\in\qri(\Omega)\Longleftrightarrow 0\in\qri(\Omega-\bar{x})=\qri(\Omega)-\{\bar{x}\},
\end{aligned}
\end{equation*}
which completes the proof of the theorem. $\h$

As we see below, in the case where $\ri(\Omega)=\emp$ the inclusions in \eqref{ri-rel1} may be strict in the simplest infinite-dimensional Hilbert space of sequences $\ell^2$ with both sets
$\iri(\Omega)$ and $\qri(\Omega)$ being nonempty; see \cite{BG,bl}.

\begin{Example}\label{ri-diff} {\rm Let $X:=\ell^2$, and let $\Omega\subset X$ be given by
\begin{equation*}
\Omega:=\Big\{x=(x_k)\in X\Big|\;\|x\|_1:=\sum_{k=1}^\infty|x_k|\le 1\Big\}.
\end{equation*}
We can  check that
\begin{equation}\label{exairi1}
\iri(\Omega)=\{x\in X\;|\;\|x\|_1<1\},
\end{equation}
and
\begin{eqnarray}\label{exaqri1}
\begin{array}{ll}
\qri(\Omega)=\Omega\setminus\big\{x=(x_k)\in X\;\big|\;&\|x\|_1=1,\\
&\exists k_0\in \N\; \mbox{\rm such that }x_k=0\;\mbox{for all}\;k\ge k_0\big\}.
\end{array}
\end{eqnarray}
To justify \eqref{exairi1}, we first take any $x\in \iri(\Omega)$ and show that $\|x\|_1<1$. Fix $u=0\in \Omega$ and find $y\in \Omega$ such that $x=ty+(1-t)u=ty$ for some $t\in (0, 1)$. Since
$\|y\|_1\leq 1$ and $0<t<1$ we have  $\|x\|_1=\|ty\|_1=t\|y\|_1<1$. To justify  the converse, we fix any $x\in X$ with $\|x\|_1<1$ and show that
\begin{equation}\label{iri l2}
\cone(\Omega-x)=\ell^1,
\end{equation}
which is a subspace of $\ell^2$. Take any $z\in \ell^1$ and choose $t>0$ sufficiently small such that $\|x\|_1+t\|z\|_1\leq 1$. It follows that $\|x+tz\|_1\leq 1$, and so $x+tz\in \Omega$. This
implies $z\in \cone(\Omega-x)$ and hence the inclusion ``$\supset$" in \eqref{iri l2}. Since the other inclusion in \eqref{iri l2} is obvious, we see that $x\in \iri(\Omega)$.

To prove \eqref{exaqri1}, observe that for any $x\in \Omega$ we have
\begin{equation}\label{nmcqri}
N(x; \Omega)=\{z\in X\; |\; \la x, z\ra =\|z\|_\infty\}.
\end{equation}
Indeed, $z\in N(x; \Omega)$ if and only if $\la z, u-x\ra \leq 0$ for all $u\in \Omega$, which is equivalent to
\begin{equation*}
\sup\{ \la z, u\ra\; |\; u\in \Omega\}=\la z, x\ra.
\end{equation*}
It is easy to check that $\Omega=\{u\in \ell^1\; |\; \|u\|_1\leq 1\}$, and so
\begin{equation*}
\sup\{ \la z, u\ra\; |\; u\in \Omega\}=\sup\{ \la z, u\ra\; |\; u\in l^1, \|u\|_1\leq 1\}=\|z\|_\infty,
\end{equation*}
which clearly implies \eqref{nmcqri}.

Now, fix any $x\in \qri(\Omega)$ and suppose to the contrary that $x$ does not belong to the set on the right-hand side of \eqref{exaqri1}. Then $x\in \ell^2$ satisfies
\begin{equation*}
\|x\|_1=1,\,x_k=0\;\mbox{for all}\;k\ge k_0\;\mbox{with some}\;k_0\in\N.
\end{equation*}
Define  $z\in \ell^2$ by $z_k=\mbox{\rm sign}(x_k)$. Then $\|z\|_\infty=1$ and $\la x, z\ra =\sum_{k=1}^\infty x_kz_k=\sum_{k=1}^\infty |x_k|=1$. It follows from \eqref{nmcqri} that $z\in N(x; \Omega)$,
which is a subspace of $X$. It follows that $-z\in N(x; \Omega)$, and so $\la -z, 0-x\ra\leq 0$. This is a contradiction since $\la -z, 0-x\ra=\la z, x\ra=1$. This contradiction tells us that $x$
belongs to the set on the right-hand side of \eqref{exaqri1}.

Next, fix any $x$ in the set on the right-hand side of
\eqref{exaqri1} and suppose again to the contrary that $x\notin
\qri(\Omega)$. Then $N(x; \Omega)$ is not a subspace of $X$. Thus,
we can find $z\neq 0$ such that $z\in N(x; \Omega)$.  It follows
from  \eqref{nmcqri} that $\la x, z\ra=\|z\|_\infty\neq 0$. Thus,
\begin{equation*} \|z\|_\infty =\la x, z\ra =\sum_{k=1}^\infty
x_kz_k\leq \sum_{k=1}^\infty |x_k| |z_k|\leq \|z\|_\infty
\sum_{k=1}^\infty |x_k|=\|z\|_\infty \|x\|_1\leq \|z\|_\infty.
\end{equation*} This implies $\|x\|_1=1$ and $|z_k|=\|z\|_\infty>0$
whenever $x_k\neq 0$. Since $z\in \ell^2$ we see that there exists
$k_0\in \N$ such that $x_k=0$ for all $k\geq k_0$. This
contradiction shows that $x\in \qri(\Omega)$ and completes the proof
of \eqref{exaqri1}.} \end{Example}

Next we present an example showing that the intrinsic relative interior may be empty for convex subsets of $\ell^2$;  see \cite{BG}.

\begin{Example}\label{iri-emp} {\rm Let $X:=\ell^2$, and let $\Omega\subset X$ be given by
\begin{equation*}
\Omega:=\Big\{x=(x_1,x_2,\ldots)\in X\;\Big|\;\|x\|_2:=\Big(\sum_{k=1}^\infty|x_k|^2\Big)^{1/2}\le 1,\,x_k\ge 0\;\mbox{for all}\;k\in\N\Big\}.
\end{equation*}
We are going to show that $\iri(\Omega)=\emp$ for this set. Assume on the contrary that there exists $\ox\in\iri(\Omega)$. Following Example \ref{ri-diff}, we see that $\|\ox\|<1$. Next, we will show
that $\ox_k>0$ for all $k\in \N$. Indeed, if, for example, $\ox_1=0$, then let $v:=(1, 0, 0, \ldots)\in X$ and easily get
\begin{equation*}
\la v, \ox\ra \leq \inf\{\la v, x\ra\; |\; x\in \Omega\}\; \mbox{\rm and } \la v, \ox\ra < \sup\{\la v, x\ra\; |\; x\in \Omega\}.
\end{equation*}
Thus, $\ox$ and $\Omega$ are properly separated, and so $\ox\in [\qri(\Omega)]^c\subset [\iri(\Omega)]^c$, a contradiction.

Proposition~\ref{iri-desc} tells us that for each $x\in\Omega$ we have $(1-\al)x+\al\ox\in\Omega$ with some $\al>1$. Fix $\ve>0$ and select an increasing sequence of natural numbers $\{k_n\}$ with
$0<\ox_{k_n}\le\ve/4^n$. Defining $\tilde x\in\Omega$ by $\tilde x_{k_n}:=\ve/2^n$ and $\tilde x_k:=0$ for all other $k\in\N$, let us check that $(1-\al)\tilde x+\al\ox\notin\Omega$ whenever $\al>1$.
Indeed, we have the estimate
\begin{equation*}
\big((1-\al)\tilde x+\al\ox\big)_{k_n}\leq \big(1-\al\big)\frac{\ve}{2^n}+\al\frac{\ve}{4^n}<0
\end{equation*}
for $n$ sufficiently large, which justifies the claim of this example.}
\end{Example}

The following result can be found in \cite[Proposition 2.16]{bl}, which gives us a condition to define a point contained in the set of quasi-relative interior points of a nonempty convex subset of an
LCTV space.

\begin{Proposition}\label{nonsupp} Let $\Omega$ be a nonempty convex subset of an LCTV space, and let $\ox\in\Omega$. Then $\ox\in\qri(\Omega)$ if and only if $\ox$ is a nonsupport point of $\Omega$.
\end{Proposition}
\noindent {\bf Proof.}  Observe first that any nonsupport point $\ox$ of $\Omega$ can be equivalently described as follows: whenever $x^*\in X^*$ we have the implication
\begin{equation}\label{nonsupp1}
\big[\la x^*,x-\ox\ra\ge 0\;\mbox{ if }\;x\in\Omega\big]\Longrightarrow\big[\la x^*,x-\ox\ra=0\;\mbox{ if }\;x\in\Omega\big].
\end{equation}
Having this in mind, assume now that $\ox\in\qri(\Omega)$. Since
\begin{equation*}
\big[\la x^*,x-\ox\ra\ge 0\;\mbox{ if }\;x\in\Omega\big]\Longrightarrow\big[\la x^*,u\ra\ge 0\;\mbox{ if }\;u\in\overline{\rm cone}(\Omega-\ox)\big]
\end{equation*}
for any $x^*\in X^*$, and since the set $\overline{\rm cone}(\Omega-\ox)$ is a subspace, we get
\begin{equation*}
\la x^*,u\ra=0\;\mbox{ if }\;u\in\overline{\rm cone}(\Omega-\ox),\;\mbox{ and so }\;\la x^*,x-\ox\ra=0\;\mbox{ if }\;x\in\Omega.
\end{equation*}
The latter means by \eqref{nonsupp1} that $\ox$ is a nonsupport point of $\Omega$.

To verify the ``if" part of the proposition, let $\ox$ be a nonsupport point of $\Omega$. Denoting $C:=\overline{\rm cone}(\Omega-\ox)$ and arguing by contradiction, suppose that $C$ is not a subspace
of
$X$, i.e., there exists $v\in C$ with $-v\notin C$. This yields by the strict separation theorem that
\begin{equation*}
\la x^*,-v\ra<\la x^*,u\ra\;\mbox{ for all }\;u\in C
\end{equation*}
for some $x^*\in X^*$. Taking into account that $C$ is a cone, we obtain that
\begin{equation*}
\la x^*,u\ra\ge 0\;\mbox{ if }\;u\in C,\;\mbox{ and so }\;\la x^*,x-\ox\ra\geq 0\;\mbox{ if }\;x\in\Omega.
\end{equation*}
On the other hand, it follows from $\la x^*,-v\ra<0$ and $v\in\overline{\rm cone}(\Omega-\ox)$ that there exist $x\in\Omega$ and $\lm > 0$ satisfying $\la x^*,-\lm(x-\ox)\ra<0$. The latter implies that
$\la x^*,x-\ox\ra>0$ contradicting~\eqref{nonsupp1}, which tells us that $\ox\in\qri(\Omega)$ and thus completes the proof of the proposition. $\h$

The following proposition is known from \cite[Lemma 3.1]{BG} and
\cite[Lemma 2.9]{bl}.  Note that an alternative proof of
this result can be obtained by using the ``convex core topology"
$\tau_c$ for which $\ri_{\tau_c}(\Omega)={\rm int}(\Omega)$; see,
e.g., \cite[Lemma~2.5]{kfgt}.

\begin{Proposition}\label{Lm_Interval} Let $\Omega$ be a convex
subset of a topological vector space $X.$ \begin{enumerate} \item If
$\ox\in\ri(\Omega)$ and $\tilde{x}\in \Omega,$ then
$(\tilde{x},\ox]\subset\ri(\Omega)$. \item If $\ox\in\iri(\Omega)$
and $\tilde{x}\in \Omega,$ then
$(\tilde{x},\ox]\subset\iri(\Omega)$. \item Suppose further that $X$
is locally convex. If $\ox\in\qri(\Omega)$ and $\tilde{x}\in
\Omega,$ then $(\tilde{x},\ox]\subset\qri(\Omega)$. \end{enumerate}
\end{Proposition} \noindent {\bf Proof.}  (a) Fix
$\ox\in\ri(\Omega)$ and $\tilde{x}\in\Omega$. If $\ox=\tilde{x}$,
then $(\tilde{x},\ox]=\{\ox\}\subset\ri(\Omega)$. Now suppose that
$\ox\ne\tilde{x}.$ Since $\ox\in\ri(\Omega)$, there exists a
neighborhood $U$ of $\ox$ such that
$U\cap\overline{\aff}(\Omega)\subset\Omega$. For any
$\lambda\in[0,1]$, we have \begin{equation*}
  (1-\lambda)\big(U\cap\overline{\aff}(\Omega)\big)+\lambda\tilde{x}\subset\Omega.
\end{equation*} Taking $y\in(\tilde{x},\ox]$, we have $y=(1-\lambda_0)\ox+\lambda_0\tilde{x}$ for some $\lambda_0\in[0,1)$. Then $V:=(1-\lambda_0)U+\lambda_0\tilde{x}$ is a neighborhood of $y$.
 Observe that
\begin{equation*}
  V\cap\overline{\aff}(\Omega)=(1-\lambda_0)\big(U\cap\overline{\aff}(\Omega)\big)+\lambda_0\tilde{x}\subset\Omega.
\end{equation*}
This implies $y\in\ri(\Omega)$, and so $(\tilde{x},\ox]\subset \ri(\Omega)$.

(b) Fix $\ox\in\iri(\Omega)$ and $\tilde{x}\in\Omega$. Without loss of generality we can assume that $\ox=0$. Fix any $\lambda\in[0,1)$ and let $x:=\lambda\tilde{x}$. In order to justify $x\in\iri(\Omega)$, we
take any $\oy\in\Omega$ with $\oy\ne x$ and show that there exists $\tilde{y}\in \Omega$ such that $x\in(\oy,\tilde{y})$. Indeed, since  $\ox=0\in\iri(\Omega)$, we have $-\alpha \oy\in\Omega$ for some $\alpha>0$. Choosing
\begin{equation*}
  \tilde{y}:=(1-\delta)\tilde{x}+\delta(-\alpha \oy)\in\Omega
\end{equation*}
where $\delta=\frac{1-\lambda}{1+\lambda\alpha}\in[0,1]$,  we can see that
\begin{equation*}
  x=(1-\gamma)\oy+\gamma \tilde{y},
\end{equation*}
where $\gamma=\frac{1+\lambda\alpha}{1+\alpha}\in(0,1)$. Thus, $x\in(\oy,\tilde{y})$ and hence Proposition~\ref{iri-desc} shows that $x\in\iri(\Omega)$.

(c) Fix $\ox\in\qri(\Omega)$ and $\tilde{x}\in\Omega$.  Let $y:=\lambda \ox+(1-\lambda)\tilde{x}$ with $\lambda\in(0,1]$. Using Proposition~\ref{qrinormal}, we will show that $N(y;\Omega)$ is a subspace.
Fix any $x^*\in N(y;\Omega)$. Then
\begin{equation}\label{eq_sec_qri}
  \la x^*,x-y\ra\leq 0\ \text{for all}\  x\in \Omega.
\end{equation} Using  $x=\ox$ in \eqref{eq_sec_qri} gives $\la x^*,\ox-\tilde{x}\ra \leq 0$. Similarly, using $x=\tilde{x}$ in \eqref{eq_sec_qri} gives $\la x^*,\tilde{x}-\ox\ra \leq 0$. From this we
conclude $\la x^*,\ox\ra =\la x^*,\tilde{x}\ra$, which together with~\eqref{eq_sec_qri} gives
\begin{equation*}
  \la x^*, x-y\ra = \la x^*, x-\tilde{x}\ra -\lambda\la x^*, \ox-\tilde{x}\ra =\la x^*,x-\ox\ra\leq 0 \ \text{for all}\  x\in \Omega.
\end{equation*}
 Thus, $x^*\in N(\ox;\Omega)$. Since $\ox\in\qri(\Omega)$, Proposition~\ref{qrinormal} shows that $-x^*\in N(\ox;\Omega)$. Since $\la x^*,\ox\ra =\la x^*,\tilde{x}\ra$, one has
\begin{equation*}
  \la -x^*,x-y\ra  = \lambda\la -x^*, x-\ox\ra+ (1-\lambda)\la -x^*,x-\tilde{x}\ra \leq  0\ \text{for all}\ x\in\Omega,
\end{equation*}
and hence $-x^*\in N(y;\Omega)$. Thus, $N(y;\Omega)$ is a subspace, and the result follows from Proposition~\ref{qrinormal}.
$\h$

The next result provides a useful version of {\em strict separation} relative to closed subspaces of Hilbert spaces; see \cite[Lemma 3.7]{CBN}.

\begin{Proposition}\label{strict-sepH} Let $L$ be a closed subspace of a Hilbert space $X$, and let $\Omega\subset L$ be a nonempty convex set with $\bar{x}\in L$ and $\bar{x}\not\in\Bar{\Omega}$.
 Then there exists  $u\in L$ such that
\begin{equation*}
\sup\big\{\langle u,x\rangle\;\big|\;x\in\Omega\big\}<\langle u,\bar{x}\rangle.
\end{equation*}
\end{Proposition}
{\bf Proof.}  Since $\ox\notin\Bar{\Omega}$, we see that the sets $\{\ox\}$ and $\Omega$ are strictly separated in $X$, which means that there exists a vector $v\in X$ such that
\begin{equation}\label{eqlm1}
\sup\big\{\la v,x\ra\;\big|\;x\in\Omega\big\}<\la v,\ox\ra.
\end{equation}
It is well known that any Hilbert space $X$ can be represented as the direct sum $X=L\oplus L^{\bot}$, where
\begin{equation*}
L^{\bot}:=\big\{w\in X\;\big|\;\la w,x\ra=0\;\text{ for all }\;x\in L\big\}.
\end{equation*}
If $v\in L^{\bot}$, then \eqref{eqlm1} immediately gives us a contradiction. Thus $v\in X$ is represented as $v=u+w$ with some $0\ne u\in L$ and $w\in L^{\bot}$. This implies that for each $x\in\Omega\subset L$ we have
\begin{equation*}
\begin{aligned}
\la u,x\ra &=\la u,x\ra+\la w,x\ra=\la v,x\ra\le\sup\big\{\la v,x\ra\;\big|\;x\in\Omega\big\}<\la v,\ox\ra\\
&=\la u+w,\ox\ra=\la u,\ox\ra, \end{aligned} \end{equation*} which
shows that $\sup\{\la u,x\ra\;|\;x\in\Omega\}<\la u,\ox\ra$. $\h$

Before establishing the next result, let us present the following useful technical lemma on intrinsic relative interiors; see \cite[Lemma 3.5]{CBN}.

\begin{Lemma}\label{lm2separate} Let $X$ be an LCTV space, and let $\Omega\subset X$ be a nonempty, closed, and convex set with $0\in\Omega\setminus\mbox{\rm iri}(\Omega)$. If ${\rm iri}(\Omega)\ne\emp$,
then $\overline{\aff}(\Omega)$ is a closed subspace of $X$ and there exists a sequence $\{x_k\}\subset-\Omega$ such that $x_k\notin\Omega$ and $x_k\to 0$ as $k\to\infty$.
\end{Lemma}
\noindent {\bf Proof.}  The set $\overline{\aff}(\Omega)$ is closed subspace of $X$ since it is a closed, affine set containing the origin. Using ${\rm iri}(\Omega)\ne\emp$ and
$0\in\Omega\setminus\mbox{\rm iri}(\Omega)$, let us show that there exists a nonzero vector $x_0\in{\rm iri}(\Omega)$ with $-tx_0\notin\Omega$ for all $t>0$. Arguing by contradiction, suppose that
 $-tx_0\in\Omega$ for some $t>0$. Then it follows from Proposition~\ref{Lm_Interval}(b) that
\begin{equation*}
0=\frac{t}{1+t}x_0+\frac{1}{1+t}\big(-tx_0\big)\in\mbox{\rm iri}(\Omega),
\end{equation*}
which clearly contradicts the assumption $0\notin{\rm iri}(\Omega)$. Letting $x_k:=-(x_0/k)\in-\Omega$ gives us that $x_k\notin{\Omega}$ for every $k$ and that $x_k\to 0$ as $k\to\infty$. $\h$

Let us now define the following property, which is automatic in finite dimensions while being very important for performing limiting procedures in infinite dimensional spaces.

\begin{Definition}\label{snc-conv} A subset $\Omega\subset X$ of a normed space $X$ is {\sc sequentially normally compact (SNC)} at $\ox\in\Omega$ if for any sequence $\{(x_k,x^*_k)\}\subset X\times X^*$
we have the implication
\begin{equation}\label{snc}
\big[x^*_k\in N(x_k;\Omega),\;x_k\in\Omega,\;x_k\to\ox,\;x^*_k\st{w^*}{\to}0\big]\Longrightarrow\|x^*_k\|\to 0.
\end{equation}
\end{Definition}

\begin{Remark}\label{snc-rem}{\rm The SNC property \eqref{snc} is  taken from \cite{m-book1} and investigated therein for general nonconvex sets in Banach spaces. In the case of closed and {\em convex}
subsets $\Omega\subset X$ of such spaces, this property can be characterized as follows \cite[Theorem~1.21]{m-book1}: If a closed and convex set $\Omega$ has nonempty relative interior, then it is SNC at
every $\ox\in\Omega$ if and only if the closure of the span of $\Omega$ is of {\em finite codimension}.}
\end{Remark}

Now we are ready to derive the aforementioned result on the quasi-regularity of convex sets in infinite dimensions; see \cite[Theorem 3.8 (d)]{CBN}.

\begin{Theorem}\label{quasi-reg} Let $\Omega\subset X$ be a nonempty, closed, and convex subset of a Hilbert space $X$. Assume in addition that $\iri(\Omega)\ne\emp$, and that $\Omega$ is SNC at every
 point $\ox\in\Omega$. Then this set is quasi-regular.
\end{Theorem}
\noindent {\bf Proof.}  First we verify that in the case where $0\notin{\rm iri}(\Omega)$ the sets $\Omega$ and $\{0\}$ are properly separated, i.e., there exists a nonzero vector $a\in X$ such that
\begin{equation}\label{pro-sep}
\sup\big\{\langle a,x\rangle\;\big|\;x\in\Omega\big\}\le 0\;\mbox{ and }\;\inf\big\{\langle a,x\rangle\;\big|\;x\in\Omega\big\}<0.
\end{equation}

If $0\not\in\Omega$, this statement is trivial. Suppose now that $0\in\Omega\setminus\text{\rm iri}(\Omega)$. Letting $L:=\overline{\mbox{\rm aff}}(\Omega)$ and employing Lemma~\ref{lm2separate} tell us
that $L$ is a subspace of $X$, and that there is a sequence $\{x_k\}\subset L$ for which $x_k\notin {\Omega}$ and $x_k\to 0$ as $k\to\infty$. By Proposition~\ref{strict-sepH} we find a sequence
$\{v_k\}\subset L$ with $v_k\ne 0$ and
\begin{equation*}
\sup\big\{\langle v_k,x\rangle\;\big|\;x\in\Omega\big\}<\la v_k,x_k\ra\;\mbox{ whenever }\;k\in\N.
\end{equation*}
Denote $w_k:=\frac{v_k}{\Vert v_k\Vert}\in L$ so $\|w_k\|=1$ for $k\in\N$ and observe that
\begin{equation}\label{eq1Lm1.5wk}
\langle w_k,x\rangle<\la w_k,x_k\ra\le\|w_k\|\cdot\|x_k\|=\varepsilon_k\;\mbox{\rm for all }\;x\in\Omega,
\end{equation}
where $\varepsilon_k:=\|x_k\|\downarrow 0$. Since $\{w_k\}$ is bounded, we let $k\to\infty$ in \eqref{eq1Lm1.5wk} and suppose without loss of generality that $w_k\xrightarrow{w}a\in L$, which yields
\begin{equation}\label{eq2Lm1.5}
\sup\big\{\langle a,x\rangle\;\big|\;x\in\Omega\big\}\le 0.
\end{equation}
To verify further the strict inequality
\begin{equation*}
\inf\big\{\langle a,x\rangle\;\big|\;x\in\Omega\big\}<0,
\end{equation*}
it suffices to show that there is $\ox\in\Omega$ with $\langle a,\ox\rangle<0$. Arguing by contradiction, suppose that $\langle a,x\rangle\ge 0$ for all $x\in\Omega$ and deduce from \eqref{eq2Lm1.5} that
$\langle a,x\rangle=0$ whenever $x\in\Omega$. Since $a\in L=\overline{\aff}(\Omega)$, there exists a sequence $a_j\to a$ as $j\to\infty$ with $a_j\in\aff(\Omega)$. The latter inclusion can be rewritten
as
\begin{equation*}
a_j=\sum_{i=1}^{m_j}\lambda^j_i\omega^j_i\;\mbox{ with }\;\sum_{i=1}^{m_j}\lambda^j_i=1\;\mbox{ and }\;\omega^j_i\in\Omega\;\mbox{ for }\;i=1,\ldots,m_j,
\end{equation*}
which readily implies the equalities
\begin{equation*}
\la a, a_j\ra=\sum_{i=1}^{m_j}\lambda^j_i\langle a,\omega^j_i\rangle=0.
\end{equation*}
The passage to the limit as $j\to\infty$ gives us $a=0$.

Next we deduce from \eqref{eq1Lm1.5wk}, by using the Br{\o}ndsted-Rockafellar theorem proved in \cite[Theorem 3.1.12]{z}, the existence of $b_k\in\Omega$ and $u_k\in X$ such that
\begin{equation}\label{eq3Lm1.5}
u_k\in N(b_k;\Omega),\;\|b_k\|\le\sqrt{\varepsilon_k},\;\mbox{ and }\;\|u_k-w_k\|\le\sqrt{\varepsilon_k}.
\end{equation}
Since $\|w_k\|=1$, it follows from \eqref{eq3Lm1.5} that $\|u_k\|\to 1$. Furthermore, we get from $w_k\xrightarrow{w}0$, $\varepsilon_k\downarrow 0$, and \eqref{eq3Lm1.5} that $u_k\xrightarrow{w}0$ as
$k\to\infty$. Remembering that $\Omega$ has the SNC property, it follows from \eqref{eq3Lm1.5} that $\|u_k\|\to 0$, which clearly contradicts the condition $\|u_k\|\to 1$ as $k\to\infty$. This tells us
that there exists $\ox\in\Omega$ such that $\la a,\ox\ra<0$. This justifies the proper separation of $\Omega$ and $\{0\}$ in \eqref{pro-sep}, so by Proposition~\ref{qri-sep1} we conclude that
$0\notin\qri(\Omega)$.

To complete the proof of the quasi-regularity of $\Omega$, it remains to show that $\qri(\Omega)\subset\iri(\Omega)$. Picking any $\ox\in\qri(\Omega)$ gives us $0\in\qri(\Omega-\ox)$. Since $\Omega$
satisfies the hypotheses of the theorem, so does $\Omega-\ox$. Applying the proof above gives $0\in\iri(\Omega-\ox)$, so $\ox\in\iri(\Omega)$. This shows that $\qri(\Omega)\subset\iri(\Omega)$ and thus finishes the proof of the theorem. $\h$

\section{Generalized Relative Interiors for Graph of Set-valued Mappings} \setcounter{equation}{0} In this section we study generalized relative interiors of  set-valued mappings with convex
graphs and give a version of Rockafellar's theorem in LCTV spaces.  Remember that a set-valued mapping $F\colon X\tto Y$ between LCTV spaces is associated with its {\em graph} \begin{equation*}
\gph(F):=\big\{(x,y)\in X\times Y\;\big|\;y\in F(x)\big\}, \end{equation*} and it is called {\em convex} if its graph is a convex subset of the product space $X\times Y$. We also consider the
{\em domain} of $F$ defined by \begin{equation*} \dom(F):=\big\{x\in X\;\big|\;F(x)\ne\emp\big\}. \end{equation*}

Let $f\colon X\to (-\infty,\infty]$ be an extended-real-valued function. Recall that the {\em domain} and the {\em epigraph} 
of  $f$ are defined by
\begin{itemize}
\item $\dom(f):=\big\{x\in X\;\big|\;f(x)<\infty\big\}$,
\item $\epi(f):=\big\{(x,\al)\in X\times\R\;\big|\;\al\ge f(x)\big\}$,
\end{itemize}
The function $f$ is said to be \emph{convex} if its epigraph is a convex set, and it is called {\em proper} if $\dom(f)\neq \emptyset$. An extended-real-valued function $g\colon X\to [-\infty,\infty)$ is
said to be {\em concave} if $-g$ is convex, and it is said to be {\em proper} if $-g$ is proper.  Define $\dom(g):=\{x\in X\; |\; -\infty<g(x)\}$. It follows from the definition that $g$ is concave if and only if  the set
\begin{equation*}
  \hypo(g):=\big\{(x,\al)\in X\times\R\;\big|\;\al\leq g(x)\big\}
\end{equation*}
is convex.

Let $f\colon X\to (-\infty,\infty]$ be an extended-real-valued function. Define the set-valued mapping $F\colon X\tto\R$ given by $F(x):=[f(x),\infty)$, then it is easy to see that $\dom(F)=\dom(f)$,
$\epi(f)=\gph(F)$, and $f$ is convex if and only if $F$ is convex.

\begin{Theorem}\label{Theo-Roc1} Let $F\colon X\tto Y$ be a convex set-valued mapping between LCTV spaces. Then we have
\begin{equation*}
  \qri(\gph(F)) \supset \{(x,y)\in X\times Y\ \big|\ x\in\qri\big(\dom(F)\big), y\in\sint \big(F(x)\big)\}.
\end{equation*}
\end{Theorem}
{\bf Proof.} Pick any $(\ox,\oy)\in X\times Y$ with $\ox\in\qri\big(\dom(F)\big)$ and $\oy\in\sint\big(F(\ox)\big)$. By a contradiction, suppose that $(\ox,\oy)\notin \qri\big(\gph(F)\big)$. Then
Proposition~\ref{qri-sep1} shows that the sets $\{(\ox,\oy)\}$ and $\gph(F)$ can be properly separated, which means that there exists $(x^*,y^*)\in X^*\times Y^*$ such that
\begin{equation}\label{Eq1-RoctheoGph}
  \la x^*,x\ra +\la y^*,y\ra \leq \la x^*,\ox\ra +\la y^*,\oy\ra \; \text{for all}\; x\in\dom(F), y\in F(x)
\end{equation}
and there exists $(\tilde{x},\tilde{y})\in \gph(F)$ such that
\begin{equation}\label{Eq2-Roctheogph}
  \la x^*,\tilde{x}\ra +\la y^*, \tilde{y}\ra < \la x^*,\ox\ra +\la y^*,\oy\ra.
\end{equation} Choosing $x=\ox$, \eqref{Eq1-RoctheoGph} shows that \begin{equation}\label{Eq3-Roctheogph} \la y^*,y\ra\leq \la y^*,\oy\ra \; \text{for all} \; y\in F(\ox). \end{equation} Since
$\oy\in\sint\big(F(\ox)\big)$, there exists a symmetric neighborhood $V$ of the origin with $V\subset \sint\big(F(\ox)\big)-\{\oy\}$. It follows from \eqref{Eq3-Roctheogph} that $\la y^*,v\ra\leq
0$ and $\la y^*,-v\ra \leq 0$ for all $0\ne v\in V$. Thus, $y^*=0$ on $V$. Since $V$ is a symmetric neighborhood of the origin, for any $y\in Y$ there exists $0\ne t\in \R$ such that $ty=v\in V$
and hence $\la y^*,y\ra =\frac{1}{t}\la y^*,v\ra=0$. Thus, $y^*=0$ on $Y$.  It follows from \eqref{Eq1-RoctheoGph} and \eqref{Eq2-Roctheogph} that the sets $\{\ox\}$ and $\dom(F)$ can be properly
separated. Proposition~\ref{qri-sep1} shows that $\ox\notin \qri\big(\dom(F)\big)$. This contradiction completes the proof of the theorem. $\h$

The following result can be found in \cite[Corollary 9(iii)]{zonuseqri}. In this paper we can see that this result is a direct corollary of Theorem~\ref{Theo-Roc1} by considering $F(x)=[f(x),\infty)$.
\begin{Corollary}\label{coro-qri-epi1}Let $f\colon X\to (-\infty,\infty]$ be a proper convex function. Then we  have
\begin{equation*}
  {\rm qri}\big(\epi(f)\big)\supset\big\{(x,\lambda)\in X\times\R\;\big|\;x\in{\rm qri}\big(\dom(f)\big),\;\lambda>f(x)\big\}.
\end{equation*}
\end{Corollary}
\begin{Theorem}\label{Theo-Roc2}
Let $F\colon X\tto Y$ be a convex set-valued mapping between LCTV spaces. Then we have
\begin{equation}\label{Eq1-Irigph}
  \iri(\gph(F)) \subset \{(x,y)\in X\times Y\ \big|\ x\in\iri\big(\dom(F)\big), y\in\iri\big(F(x)\big)\}.
\end{equation} \end{Theorem} {\bf Proof.} To prove \eqref{Eq1-Irigph},  pick any $(\ox,\oy)\in\iri\big(\gph(F)\big)$, we first check that $\ox\in\iri\big(\dom(F)\big)$. For any $x\in\dom(F)$ with
$x\ne\ox$ and $y\in F(x),$ one has $(\ox,\oy)\ne
(x,y)\in\gph(F)$. Since $(\ox,\oy)\in\iri\big(\gph(F)\big)$, Proposition~\ref{iri-desc} shows that there exists $(u,v)\in\gph(F)$ such that
$(\ox,\oy)\in\big((x,y),(u,v)\big)$ and hence $\ox\in(x,u)$. Applying Proposition~\ref{iri-desc} gives us that $\ox\in\iri\big(\dom(F)\big)$.

Let us now show that $\oy\in\iri\big(F(\ox)\big)$. Pick any $\tilde{y}\in F(\ox)$ with $\tilde{y}\ne\oy$, and hence $(\ox,\oy)\ne(\ox,\tilde{y})$. By the assumption that
$(\ox,\oy)\in\iri\big(\gph(F)\big)$, Proposition~\ref{iri-desc} shows that there exists $(u,v)\in\gph(F)$ such that $(\ox,\oy)\in\big((\ox,\tilde{y}),(u,v)\big)$. Thus, $\ox\in (\ox,u)$ and hence
$u=\ox$. Therefore, $v\in F(u)=F(\ox)$. This means that there exists $v\in F(\ox)$ such that $\oy\in (\tilde{y},v)$. Applying Proposition~\ref{iri-desc} gives us that $\oy\in\iri\big(F(\ox)\big)$ and
\eqref{Eq1-Irigph} was proved.
$\h$
\begin{Theorem}\label{Theo-Roc3} Let $F\colon X\tto Y$ be a convex set-valued mapping between LCTV spaces. Suppose that $\gph(F)$ is quasi-regular and $\sint\big(F(x)\big)\ne \emptyset$ for all
$x\in \dom (F)$. Then we have
\begin{equation}\label{Eq1-The-Roc3}
  \qri(\gph(F)) = \{(x,y)\in X\times Y\ \big|\ x\in\qri\big(\dom(F)\big), y\in\sint \big(F(x)\big)\}
\end{equation} and $\dom (F)$ is quasi-regular. \end{Theorem} {\bf
Proof.} Set \begin{equation*} \Omega:=\{(x,y)\in X\times Y\ \big|\
x\in\qri\big(\dom(F)\big), y\in\sint \big(F(x)\big)\}.
\end{equation*} Applying Theorem~\ref{Theo-Roc1}, we have
\begin{equation}\label{Eq2-Theo-Roc3}
 \qri\big(\gph(F)\big)\supset \Omega.
\end{equation}
Since $\gph(F)$ is quasi-regular and $\sint\big(F(x)\big)\ne\emptyset$  for all $x\in\dom(F)$, Theorem~\ref{ri-rel} and Theorem~\ref{Theo-Roc2} show that
\begin{equation}\label{Eq3-Theo-Roc3}
  \qri\big(\gph(F)\big)=\iri\big(\gph(F)\big)\subset \{(x,y)\in X\times Y\ \big|\ x\in\iri\big(\dom(F)\big), y\in\iri\big(F(x)\big)\}\subset \Omega.
\end{equation}
The inclusions  \eqref{Eq2-Theo-Roc3} and \eqref{Eq3-Theo-Roc3} imply \eqref{Eq1-The-Roc3}. It also follows from \eqref{Eq2-Theo-Roc3} and \eqref{Eq3-Theo-Roc3} that
$\iri\big(\dom(F)\big)=\qri\big(\dom(F)\big)$ and hence $\dom(F)$ is quasi-regular.
$\h$

 The following result is a direct corollary of Theorem~\ref{Theo-Roc3} by considering $F(x)=[f(x),\infty)$. This result also can be seen in \cite[Theorem 5.6]{CBN}.
 \begin{Corollary}\label{coro-qri-epi2} Let $f\colon X\to (-\infty,\infty]$ be a proper convex function. If  $\epi(f)$ is quasi-regular, then $\dom(f)$ is quasi-regular and
\begin{equation*}\label{qriepi}
{\rm qri}\big(\epi(f)\big)=\big\{(x,\lambda)\in X\times\R\;\big|\;x\in{\rm qri}\big(\dom(f)\big),\;\lambda>f(x)\big\}.
\end{equation*}
 \end{Corollary}

In the following preliminary material leading to Theorem~\ref{Theo-Iri-C_epi}, we use terminology and apply results from~\cite{luxemburg-1971-riesz}. Let $Y$ be a linear space and $C$ be a nonempty
convex cone in $Y$. We define a relation in $Y$ by
\begin{equation*}
  y_1\leq_{C}y_2\; \text{if and only if}\; y_2-y_1\in C.
\end{equation*}
 This defines a partial ordering on $Y$, and in this case we say that $Y$ is partially ordered by the nonempty convex cone $C$ and refer to $Y$ as an ordered vector space. We also define the relation
\begin{equation*}
  y_1<_{C}y_2\; \text{if and only if}\; y_2-y_1\in C\; \text{and} \; y_1\ne y_2.
\end{equation*}
We say that $Y$ is {\em totally ordered} by the nonempty convex cone $C$ if for any $y_1,y_2\in Y$ we have either $y_1\leq_Cy_2$ or $y_2\leq_Cy_1.$

It is immediate that $0\leq_C y$ if and only if $y\in C$, so we refer to the points in $C$ as the {\em nonnegative elements} and use the notation $Y^+$. We refer to the supremum and infimum of a subset
$U\subset Y$ in the usual sense, namely $\bar y=\sup U$ if $u\leq_C \bar y$ holds for all $u\in U$ and $\bar y\leq_C y$ for any other upper bound $y$ for $U$, and similarly for $\inf U$.
If $\{y_n\}_{n\in\N}$ is a sequence in $Y$, we naturally refer to it is a decreasing sequence if $y_{n+1}\leq_C y_n$ for every $n\in\N$. Finally, a decreasing, nonnegative sequence $\{p_n\}_{n\in\N}\subset Y^+$ {\em converges in order} to 0 if $0=\inf_{n\in\N}\{p_n\}$ in which case we write $p_n\downarrow_C 0$. Such sequences are used to define a general notion of order convergence $y_n\rightarrow_C y$ in ordered vector spaces~\cite[Theorem~16.1]{luxemburg-1971-riesz}.

Let $f\colon X\to Y$ be a function between LCTV spaces, where $Y$ is partially ordered by a nonempty convex cone $C$. We define the $C$-{\em domain}  and the $C-${\em epigraph}  of $f$ as follows:
\begin{equation*}
  \dom_C(f):=\{x \in X \; \big| \; \exists y\in Y \; \text{with} \; f(x)<_Cy\},
\end{equation*}
\begin{equation*}
  \epi_C(f):=\{(x,y)\in X\times Y\; \big| \; x\in\dom_C(f),f(x)\leq_{C}y\}.
\end{equation*} The function $f$ is said to be $C$-{\em proper} if
$\dom_C(f)\ne\emptyset$. We say that $f$ is $C$-convex if for all
$x_1,x_2\in\dom_C(f)$ and $\lambda\in[0,1]$ one has
\begin{equation*} f\big(\lambda x_1+(1-\lambda)x_2\big)\leq_C
\lambda f(x_1)+(1-\lambda)f(x_2). \end{equation*}  The
next lemma is taken from \cite[Proposition 6.2]{dtl}.

\begin{Lemma} Let $X$ and $Y$ be LCTV spaces with $Y$ is partially
ordered by a nonempty convex cone $C$ and let $f\colon X\to Y$ be a
$C$-proper function.  Then $f$ is $C$-convex if and only if
 $\epi_C(f)$ is a convex subset of $X\times Y$.
\end{Lemma}

The following result gives us a representation of the intrinsic interior of $C$-epigraph for a $C$-proper convex function between LCTV spaces. To achieve equality, we restrict our attention to
{\em Archimedean} ordered vector spaces which have the property that $n^{-1}y\downarrow_C 0$ whenever $y\in Y^+$. This is not true in general even if $Y$ is totally ordered. For example, consider
$Y=\R^2$ with the lexicographical partial ordering which is induced by the nonempty convex cone
\begin{equation}\label{eq:lexi_cone}
C=\{(x,y)\in\R^2\;\big|\;0<x\}\cup\{(0,y)\in\R^2\;\big|\;0\leq y\}.
\end{equation}
In this case, the sequence corresponding to $p_n=n^{-1}(1,1)$ where $n\in\N$ is a decreasing, nonnegative sequence but $(0,0)<_C(0,m)\leq_C p_n$ for every $m,n\in\N$. The origin $(0,0)$ is a
lower bound for for $\{p_n\}$ but $(0,0)\not=\inf_{n\in\N}\{p_n\}$. In fact, $\inf_{n\in\N}\{p_n\}$ does not exist even though the sequence $\{p_n\}$ is bounded below. This example shows that in a
non-Archimedean ordered vector space $Y$, it is not true in general that $y_n\rightarrow_C y$ in $Y$ and $\alpha_n\rightarrow\alpha$ in $\R$ gives $\alpha_n y_n\rightarrow_C\alpha y$ in $Y$. However, this {\em is} true if $Y$ is Archimedean~\cite[Section~16]{luxemburg-1971-riesz}, and we will use this property in the next result.

\begin{Theorem}\label{Theo-Iri-C_epi} Let $X$ and $Y$ be LCTV spaces, let $Y$ be partially ordered by a nonempty convex cone $C$, and let $f\colon X\to Y$ be a $C$-proper convex function. Then we have
\begin{equation}\label{priepi}
  \iri\big(\epi_C(f)\big)\subset\big\{(x,y)\in X\times Y\;\big|\;x\in \iri\big(\dom_C(f)\big),\;f(x)<_C y\big\}.
\end{equation}
If in addition $Y$ is totally ordered by $C$ and is Archimedean, then
\begin{equation}\label{priepi2}
  \iri\big(\epi_C(f)\big)=\big\{(x,y)\in X\times Y\;\big|\;x\in \iri\big(\dom_C(f)\big),\;f(x)<_C y\big\}.
\end{equation}
\end{Theorem}
{\bf Proof.}
Denoting by $\Omega$ the set on the right-hand side of \eqref{priepi}. Let us first verify the inclusion \eqref{priepi}. Pick any $(\ox,\oy)\in{\rm iri}(\epi_C(f))$ and first check that
$\ox\in{\rm iri}\big(\dom_C(f)\big)$. Fixing $x\in\dom_C(f)$ with $x\ne\ox$, we get $(x,y)\in\epi_C(f)$, where $y:=f(x)$. Then Proposition~\ref{iri-desc} ensures the existence of $(u,v)\in\epi_C(f)$
such that
\begin{equation*}
(\ox,\oy)\in\big((x,y),(u,v)\big),
\end{equation*}
which shows that $\ox\in(x,u)$. Applying Proposition~\ref{iri-desc} again yields $\ox\in{\rm iri}\big(\dom_C(f)\big)$.

Let us now show that $f(\ox)<_C\bar{y}$. Arguing by contradiction, suppose that $\oy=f(\ox)$ and take any $(\ox,\tilde{y})\in\epi_C(f)$ with $f(\ox)<_C\tilde{y}$. Thus
$(\ox,\tilde{y})\ne(\ox,\bar{y})=\big(\ox,f(\ox)\big)$. Then it follows from Proposition~\ref{iri-desc} that there exists $(\bar{u},\bar{v})\in\epi_C(f)$ such that
$\big(\ox,f(\ox)\big)\in\big((\ox,\tilde{y}),(\bar{u},\bar{v})\big)$, and hence we can find  $t_0\in(0,1)$ such that \begin{equation*} \ox=t_0\ox+(1-t_0)\bar{u}\;\mbox{ and
}\;\bar{y}=t_0\tilde{y}+(1-t_0)\bar{v}. \end{equation*} Employing the $C$-convexity of $f$ shows that \begin{equation*} t_0\tilde{y}+(1-t_0)\bar{v}=\bar{y}=f(\ox)\le_C
t_0f(\ox)+(1-t_0)f(\bar{u})<_Ct_0\tilde{y}+(1-t_0)f(\bar{u}) \end{equation*} thus verifying that $\bar{v}<_Cf(\bar{u})$, and hence $(\bar{u},\bar{v})\notin\epi_C(f)$. The obtained contradiction
tells us that $f(\ox)<_C\bar{y}$ and therefore justifies \eqref{priepi}.

We now suppose that $Y$ is totally ordered by $C$ and is Archimedean. To prove~\eqref{priepi2}, fix any $(\ox,\bar{y})\in\Omega$ giving us $\ox\in{\rm iri}(\dom_C(f))$ and $f(\ox)<_C\bar{y}$.
Picking now any $(x,y)\in\epi_C(f)$ with $(x,y)\ne(\ox,\bar{y})$, let us verify the existence of $(\bar{u},\bar{v})\in\epi_C(f)$ for which
\begin{equation*}
(\ox,\bar{y})\in\big((x,y),(\bar{u},\bar{v})\big).
\end{equation*}
To proceed, we consider following two cases:\\[1ex]
{\bf Case 1: $x\ne\ox$}. Since $\ox\in{\rm iri}(\dom_C(f))$ and $\ox\ne x\in\dom_C(f)$, there exists $\tilde{x}\in\dom_C(f)$ such that $\ox\in(x,\tilde{x})$. Choose $\tilde{y}\in Y$ satisfying
\begin{equation}\label{eq:xbar_ybar}
(\ox,\bar{y})\in\big((x,y),(\tilde{x},\tilde{y} )\big)
\end{equation}
and we will find $(\bar{u},\bar{v})\in((\ox,\bar{y}),(\tilde{x},\tilde{y}))$ with $(\bar{u},\bar{v})\in\epi_C(f)$. Arguing by contradiction, suppose that for every
$(u,v)\in((\ox,\bar{y}),(\tilde{x},\tilde{y} ))$ we have $(u,v)\notin\epi_C(f)$, i.e., $v<_Cf(u)$ since $Y$ is totally ordered. Fix any $t\in(0,1)$ and define the $t$-dependent elements
\begin{equation*}
u_t:=t\tilde{x}+(1-t)\ox\;\mbox{ and }\;v_t:=t\tilde{y}+(1-t)\bar{y},
\end{equation*}
for which we get $(u_t,v_t)\in((\ox,\bar{y}),(\tilde{x},\tilde{y}))$. The $C$-convexity of $f$ ensures that
\begin{equation*}
t\tilde{y}+(1-t)\bar{y}=v_t<_Cf(u_t)\le_C tf(\tilde{x})+(1-t)f(\ox)\le_C t f(\tilde{x})+(1-t)\bar{y}.
\end{equation*}
Since $Y$ is Archimedean, letting $t\downarrow 0$ shows that $\bar{y}=f(\ox)$, a contradiction that verifies the existence of a pair $(\bar{u},\bar{v})$ with
\begin{equation}\label{eq:ubar_vbar}
(\bar{u},\bar{v})\in\big((\ox,\bar{y}),(\tilde{x},\tilde{y})\big)
\end{equation}
and $(\bar{u},\bar{v})\in\epi_C(f)$. Equations~\eqref{eq:xbar_ybar} and~\eqref{eq:ubar_vbar} give that $(\ox,\bar{y})\in((x,y),(\bar{u},\bar{v}))$, and so it follows from Proposition~\ref{iri-desc}
that $(\ox,\bar{y})\in{\rm iri}(\epi_C(f))$. Therefore, we have \eqref{priepi2}.\\[1ex]
{\bf Case~2: $x=\ox$}. Since $(x,y)\ne (\ox,\oy)$, we have $y\ne \oy$. We will show that there exists $\bar v\in Y$ such that $(\ox,\bar v)\in\epi_C(f)$ and $\oy\in (y,\bar v)$. To arrive at a
contradiction, suppose that for any $v\in Y$ with $\oy\in(y,v)$ we have $(\ox,v)\notin \epi_C(f)$, i.e. $v<_Cf(\ox)$ since $Y$ is totally ordered. Fix any $t\in(0,1)$ and define the $t$-dependent element
\begin{equation*}
  v_t=\frac{1}{t}\oy +\frac{t-1}{t}y.
\end{equation*} Then $\oy\in (y,v_t)$ and hence \begin{equation*} \frac{1}{t}\oy +\frac{t-1}{t}y =v_t<_Cf(\ox). \end{equation*} Again since $Y$ is Archimedean, letting $t\to1$ shows that
$\oy\leq_C f(\ox)$. This contradiction shows that there exists $(\ox,\bar{v})\in \epi_C(f)$ such that $(\ox,\oy)\in\big((x,y),(\ox,\bar{v})\big)$ and hence $(\ox,\oy)\in\iri\big(\epi_C(f)\big)$
and the theorem is proved. $\h$

The following counterexamples show that the equality~\eqref{priepi2} can fail if an ordered vector space is not totally ordered or is not Archimedean.
\begin{Example}\rm Consider $f\colon \R\to\R^2$ given by $f(x):=(0,0)$ for all $x\in \R$ and the nonempty convex cone
 \begin{equation*}
  C:=\R^2_+=\{(y,z)\in \R^2 \; \big| \; 0\leq y,0\leq z\}.
\end{equation*} Suppose that $\R^2$ is partially ordered by $C$. Then we can see that $\dom_C(f)=\R$, $\epi_C(f)=\R\times C$, and the resulting ordering is not total however $\R^2$ with this
partial ordering is Archimedean. One has \begin{equation*}
  \begin{aligned}
 \iri\big(\epi_C(f)\big)&=\R\times \sint C=\R\times\{(y,z)\in\R^2 \; \big| \; 0<y,0<z\}\\
 &\ne \R\times\left(C\setminus\{(0,0)\}\right)\\
 &= \R\times\{(y,z)\in\R^2 \; \big| \; 0\leq y,0\leq z, (y,z)\ne(0,0)\}\\
 &= \{(x,y,z)\in\R^3 \; \big| \; x\in\iri\big(\dom_C(f)\big),f(x)<_C(y,z)\}.
  \end{aligned}
\end{equation*}

Furthermore, consider the same function $f$ and the lexicographical partial ordering on $\R^2$, which is induced by the nonempty convex cone $C$ defined in~\eqref{eq:lexi_cone}. Here we have
$\dom_C(f)=\R$, $\epi_C(f)=\R\times C$, and the resulting ordered vector space is totally ordered while not being Archimedean. In this case, we have \begin{equation*}
  \begin{aligned}
 \iri\big(\epi_C(f)\big)&=\R\times \sint C=\R\times\{(y,z)\in\R^2 \; \big| \; 0<y\}\\
 &\ne \R\times\left(C\setminus\{(0,0)\}\right)\\
 &= \R\times\big(\{(x,y)\in\R^2\;\big|\;0<x\}\cup\{(0,y)\in\R^2\;\big|\;0<y\}\big)\\
 &= \{(x,y,z)\in\R^3 \; \big| \; x\in\iri\big(\dom_C(f)\big),f(x)<_C(y,z)\}.
  \end{aligned}
\end{equation*}

\end{Example}

The following result is a direct corollary of Theorem~\ref{Theo-Iri-C_epi} by considering $Y=\R$, $C=[0,\infty)$. This result also can be seen in \cite[Theorem 5.4]{CBN}.

\begin{Corollary}\label{iri-qri-epi} Let $f\colon X\to (-\infty,\infty]$ be a proper convex function. Then we  have
\begin{equation*}
  \iri\big(\epi(f)\big)=\big\{(x,\lambda)\in X\times\R\;\big|\;x\in \iri\big(\dom(f)\big),\;\lambda>f(x)\big\}.
\end{equation*}
\end{Corollary}

\section{ Convex Separation via Extended Relative Interiors and Fenchel-Rockafellar Theorem in LCTV Spaces}
\setcounter{equation}{0}

In this section we derive enhanced versions of convex separation theorems for nonsolid sets in LCTV spaces under extended relative interiority assumptions. Then we study the Fenchel-Rockafellar theorem in
the general framework of LCTV spaces. Our goal is to
  establish new sufficient conditions which guarantee the validity of \eqref{Eq-Fen-Roc} for functions defined on LCTV spaces.

Before establishing the main convex separation theorem for two nonsolid sets in terms of their extended relative interiors in LCTV spaces, we present some {\em calculus rules} involving both intrinsic
relative and quasi-relative interiors of convex sets. These results were proved
for Banach spaces in \cite[Lemma 3.3, Lemma 3.4, and Lemma 3.6]{BG}. These rules are of their own interest while being instrumental to derive the aforementioned convex separation theorem.
\begin{Theorem}\label{qri-calculus} Let $A\colon X\to Y$ be a linear continuous mapping between two LCTV spaces, and let $\Omega\subset X$ be a convex set. The following assertions hold:
\begin{enumerate}
\item  $A\big(\iri(\Omega)\big)\subset\iri\big(A(\Omega)\big)$;  moreover, equality holds if $\iri(\Omega)\ne\emp$.
\item $A\big(\qri(\Omega)\big)\subset\qri\big(A(\Omega)\big)$; moreover, equality holds if $\qri(\Omega)\ne\emp$ and the set $A(\Omega)$ is quasi-regular.
\end{enumerate}
\end{Theorem}
\noindent {\bf Proof.} It follows from Proposition \ref{iri-desc} that $\iri(\Omega)=\icr(\Omega)$, where $\icr(\Omega)$ is the relative algebraic interior of $\Omega$. Therefore, the proof of part (a) can be found in \cite[Proposition 6.3.2]{ktz} or \cite[Proposition 2.1 and Corollary 2.1]{NZ}.  To verify assertion (b) pick any $\ox\in\qri(\Omega)$ and deduce from \eqref{qri-char} that
$N(\ox;\Omega)$ is a subspace of $X^*$. Then take $y^*\in N(A(\ox);A(\Omega))$ meaning that $\la y^*,A(x)-A(\ox)\ra\le 0$ for all $x\in\Omega$, which tells us that $A^*(y^*)\in N(\ox;\Omega)$. By the
subspace property of $N(\ox;\Omega)$ we get that $-A^*(y^*)\in N(\ox;\Omega)$, which is equivalent to $-y^*\in N(A(\ox);A(\Omega))$. Thus the normal cone $N(A(\ox);A(\Omega))$ is a subspace of $Y^*$, and
so $A(\ox)\in\qri(A(\Omega))$ by \eqref{qri-char}. This verifies the inclusion ``$\subset$" in (b).

Finally, we prove the inclusion ``$\supset$" in assertion (b) under the assumptions that $\qri(\Omega)\ne\emptyset$ and $A(\Omega)$ is quasi-regular. Fix  $\ox\in\qri(\Omega)$ and set $\bar{y}:=A(\ox)$.
By the forward inclusion in~(b) which we proved above, we have  $\bar y=A(\ox)\in A\big(\qri(\Omega)\big) \subset  \qri\big(A(\Omega)\big)$. Fix any $y\in \qri\big(A(\Omega)\big)=\iri\big(A(\Omega)\big)$. If $y=\bar{y},$
then $y\in A\big(\qri(\Omega)\big)$. If $y\ne\bar{y},$ by Proposition \ref{iri-desc}  there exists $u\in A(\Omega)$ such that $y\in (u,\bar{y})$. Pick $\tilde{x}\in \Omega$ such that $A(\tilde{x})=u$ and get
\begin{equation*}
y=tu+(1-t)\bar{y}=tA(\tilde{x})+(1-t)A(\ox)=A(t\tilde{x}+(1-t)\ox)
\end{equation*}
for some $t\in (0, 1)$.   Then $y=A(x_t)$, where $x_t:=t\tilde{x}+(1-t)\ox\in\qri(\Omega)$  by Proposition~\ref{Lm_Interval}(c). Thus, $y\in A\big(\qri(\Omega)\big)$, which completes the proof. $\h$

The next theorem presents the major separation result for two nonsolid convex sets in arbitrary LCTV spaces.

\begin{Theorem}\label{Theoqricap} Let $\Omega_1$ and $\Omega_2$ be
convex subsets of an LCTV space $X$ such that $\Omega_1\cap \Omega_2\neq\emptyset$. Assume that
$\qri(\Omega_1)\ne\emp$, $\qri(\Omega_2)\ne\emp$, and the set
difference $\Omega_1-\Omega_2$ is quasi-regular. Then the sets
$\Omega_1$ and $\Omega_2$ are properly separated if and only if
\begin{equation}\label{qri-sep}
\qri(\Omega_1)\cap\qri(\Omega_2)=\emp. \end{equation} \end{Theorem}
{\bf Proof.}  First we verify that the assumptions of the theorem
ensure that \begin{equation}\label{qri-diff}
\qri(\Omega_1-\Omega_2)=\qri(\Omega_1)-\qri(\Omega_2).
\end{equation} Indeed, define a linear continuous mapping $A\colon
X\times X\to X$ by $A(x,y):=x-y$ and let
$\Omega:=\Omega_1\times\Omega_2$. It is easy to check that
$\qri(\Omega)=\qri(\Omega_1)\times\qri(\Omega_2$), and thus
$\qri(\Omega)\ne\emp$ under the assumptions made. Applying the
equality condition from Theorem~\ref{qri-calculus}(b) yields
\begin{equation*} \mbox{\rm
qri}(\Omega_1-\Omega_2)=\qri\big(A(\Omega)\big)=A\big(\qri(\Omega)\big)=\qri(\Omega_1)-\mbox{\rm
qri}(\Omega_2), \end{equation*} and thus we arrive at the claimed
equality \eqref{qri-diff}.

Consider further the set difference $\Omega:=\Omega_1-\Omega_2$ and get from \eqref{qri-diff} that condition \eqref{qri-sep} reduces to
\begin{equation*}
0\notin\qri(\Omega_1-\Omega_2)=\qri(\Omega_1)-\qri(\Omega_2),
\end{equation*}
and hence $0\notin\qri(\Omega_1-\Omega_2)=\qri(\Omega)$ under the fulfillment of \eqref{qri-sep}. Since $0\in \Omega=\Omega_1-\Omega_2$ due to the assumption $\Omega_1\cap\Omega_2\neq\emptyset$,  Proposition~\ref{qri-sep1} tells us that the sets $\Omega$ and $\{0\}$ are properly separated,
which clearly ensures the proper separation of the sets $\Omega_1$ and $\Omega_2$.

To verify the opposite implication, suppose that $\Omega_1$ and $\Omega_2$ are properly separated, which implies that the sets $\Omega=\Omega_1-\Omega_2$ and $\{0\}$ are properly separated as well.
Then using Proposition~\ref{qri-sep1} and Theorem~\ref{qri-calculus} yields
\begin{equation*}
0\notin\qri(\Omega)=\qri(\Omega_1-\Omega_2)=\qri(\Omega_1)-\qri(\Omega_2),
\end{equation*}
and thus $\qri(\Omega_1)\cap\qri(\Omega_2)=\emp$, which completes the proof. $\h$

The following result presents a proper separation theorem in LCTV
spaces via relative interior.

\begin{Corollary}\label{Theo_Ri_separare} Let $\Omega_1$ and
$\Omega_2$ be convex subsets of an LCTV space $X$ such that $\Omega_1\cap\Omega_2\neq\emptyset$. Suppose that
$\ri(\Omega_1)\ne\emptyset$,  $\ri(\Omega_2)\ne\emptyset$, and
$\ri(\Omega_1-\Omega_2)\ne\emptyset$. Then $\Omega_1$ and $\Omega_2$
can  be properly separated if and only if \begin{equation*}
\ri(\Omega_1)\cap\ri(\Omega_2)=\emptyset. \end{equation*}
\end{Corollary} {\bf Proof}. Applying Theorem~\ref{ri-rel}, we have
that $\qri(\Omega_1)=\iri(\Omega_1)=\ri(\Omega_1)\ne\emptyset$,
$\qri(\Omega_2)=\iri(\Omega_2)=\ri(\Omega_2)\ne\emptyset$, and the
set $\Omega_1-\Omega_2$ is quasi-regular. Applying
Theorem~\ref{Theoqricap} gives the conclusion of this
corollary. $\h$

 Note that Theorem~6.2.3 from \cite{ktz} establishes
sufficient conditions for the proper separation of convex sets in
LCTV spaces assuming that $\iri(\Omega_1)\ne\emp$,
$\iri(\Omega_2)\ne\emp$, and
$\iri(\Omega_1)\cap\iri(\Omega_2)=\emp$.

As seen above and will be seen in the sequel, the
quasi-regularity of convex sets is needed for the fulfillment of
many important results. Theorem~\ref{ri-rel} tells us the
quasi-regularity of a convex set $\Omega$ holds in LCTV spaces if
$\ri(\Omega)\ne\emp$ (in particular, for nonempty convex sets in
finite-dimensions), and of course if $\Omega$ is a solid convex set.

Now we recall the definition of the Fenchel conjugate of extended-real-valued functions. This concept plays a central role in convex analysis and convex optimization.
\begin{Definition} Let $f\colon X\to (-\infty,\infty]$ be a convex function and  let $g\colon X\to [-\infty,\infty)$ be a concave function.
\begin{enumerate}
  \item The {\sc convex Fenchel conjugate} of $f$ is  the function $f^*\colon X^*\to[-\infty,\infty]$ given by
\begin{equation*}
  f^*(x^*):=\sup\{\la x^*,x\ra -f(x)\ |\ x\in X\}, \; x^*\in X^*.
\end{equation*}
  \item The {\sc concave Fenchel conjugate} of $g$ is the  function $g_*\colon X^*\to[-\infty,\infty]$ given by
\begin{equation*}
  g_*(x^*):=\inf\{\la x^*,x\ra -g(x)\ |\ x\in X\}, \; x^*\in X^*.
\end{equation*}
\end{enumerate}
\end{Definition}
Note that if $\dom (f)\ne \emptyset,$ then $f^*\colon X^*\to (-\infty,\infty]$ is a convex function and
\begin{equation*}
  f^*(x^*)=\sup\{\la x^*,x\ra -f(x)\ |\ x\in\dom(f)\},\; x^*\in X^*.
\end{equation*}

Before the formulation and proof of the main duality theorem given below, we present the following simple lemma about some properties of intrinsic relative and quasi-relative interiors as well as
quasi-regularity of convex sets that are taken from Definition~\ref{qri}.

\begin{Lemma}\label{Qri_trans} Let $\Omega$ be a convex subset of an LCTV space $X$, and let $q\in X$. Then we have: \begin{enumerate} \item $\iri(q+\Omega)=q+\iri(\Omega)$. \item
$\qri(q+\Omega)=q+\qri(\Omega).$ \item $\Omega$ is quasi-regular if and only if $\Omega+q$ is quasi-regular. \end{enumerate} \end{Lemma}
{\bf Proof.} Fix any $x\in\Omega$ and observe easily that
\begin{equation*} \cone(q+\Omega-x)=\cone\big(\Omega-(x-q)\big)\;\mbox{and}\;\overline{\cone}(q+\Omega-x)=\overline{\cone}\big(\Omega-(x-q)\big). \end{equation*} Then we deduce from the
definitions of iri and qri that $x\in\iri(q+\Omega)$ if and only if $x-q\in\iri(\Omega)$, and that $x\in\qri(q+\Omega)$ if and only if $x-q\in\qri(\Omega)$. This readily verifies both assertions
(a) and (b). Assertion (c) follows directly from (a) and (b) and the definition of quasi-regularity. $\h$

It is more convenient in what follows to consider the primal
optimization problem in the following {\em difference form}:
\begin{equation}\label{FPP-diff} \mbox{minimize }\;f(x)-g(x)\;\mbox{
subject to }\;x\in X, \end{equation} where $f\colon
X\to(-\infty,\infty]$ is a {\em proper convex} function, while
$g\colon X\to[-\infty,\infty)$ is a {\em proper concave} function,
i.e., such that the function $-g$ is a proper convex function.  Note that
\eqref{FPP-diff} is a minimization problem with the convex objective
$f+(-g)$.

Now we are ready to establish the aforementioned duality theorem for problem \eqref{FPP-diff} written in the difference form.

\begin{Theorem}\label{TheoFenchel} Let $f\colon X\to(-\infty,\infty]$ be a proper convex function, and let $g\colon X\to [-\infty,\infty)$ be a proper concave function defined on an LCTV space $X$ such that  exists $\hat{x}\in X$ satisfying $\inf\big\{f(x)-g(x)\;\big|\;x\in X\big\}=f(\hat{x})-g(\hat{x})$.
Then we have the duality relationship \begin{equation}\label{duality_Theo1} \inf\big\{f(x)-g(x)\;\big|\;x\in X\big\}=\sup\big\{g_*(x^*)-f^*(x^*)\;\big|\;x^*\in X^*\big\} \end{equation} provided
that the following conditions are satisfied  simultaneously: \begin{enumerate} \item $\qri\big(\dom(f)\big)\cap\qri\big(\dom(g)\big)\ne\emp$. \item All the convex sets $\dom(f)-\dom(g)$,\;
$\epi(f)$, and $\epi(f)-\hypo(g)$ are quasi-regular.
\end{enumerate} \end{Theorem}
{\bf Proof.} Observe first that for any $x\in X$ and $x^*\in X^*$ we have the inequalities \begin{equation*}
f(x)+f^*(x^*)\ge\la x^*,x\ra\ge g(x)+g_*(x^*), \end{equation*} which immediately yield the estimate \begin{equation*} \inf\big\{f(x)-g(x)\;\big|\;x\in
X\big\}\ge\sup\big\{g_*(x^*)-f^*(x^*)\;\big|\;x^*\in X^*\big\}. \end{equation*} Denoting $\alpha:=\inf\{f(x)-g(x)\ |\ x\in X\}$, it is easy to see that \eqref{duality_Theo1} holds if
$\alpha=-\infty$. Considering the case where $\al$ is finite, we are going to show that there exists $\ox^*\in X^*$ such that $g_*(\ox^*)-f^*(\ox^*)\ge\alpha$, which would readily justify
\eqref{duality_Theo1}. To proceed, define the sets \begin{equation*} \Omega_1:=\epi(f)\;\mbox{ and }\;\Omega_2:=\big\{(x,\mu)\in X\times\R\;\big|\;\mu\le g(x)+\alpha\big\}. \end{equation*} 
It can be verified that $\Omega_1\cap\Omega_2\neq\emptyset$. Since
the set $\epi(f)$ is quasi-regular, we get by Corollary~\ref{coro-qri-epi1} and Corollary~\ref{coro-qri-epi2} that
\begin{equation*} \qri(\Omega_2)\supset\big\{(x,\mu)\in X\times \R\ \big|\ x\in\qri\big(\dom(g)\big),\ \mu< g(x)+\alpha\big\}, \end{equation*}
\begin{equation*} \qri(\Omega_1)=\big\{(x,\lambda)\in X\times \R\ \big|\ x\in\qri\big(\dom(f)\big),\ f(x)<
\lambda\big\}. \end{equation*}  It follows from the
qualification condition $\qri\big(\dom(f)\big)\cap\qri\big(\dom(g)\big)\ne\emp$ in (a) that $\qri(\Omega_1)\ne\emp$ and $\qri(\Omega_2)\ne\emp$. Using $f(x)\ge g(x)+\alpha$ for all $x\in X$ yields \begin{equation*}
\qri(\Omega_1)\cap\Omega_2=\emp,\;\mbox{ and so }\;\qri(\Omega_1)\cap\qri(\Omega_2)=\emp. \end{equation*} Observing further that $\Omega_2=\hypo(g)+\{(0,\alpha)\}$, we get
\begin{equation*}\Omega_1-\Omega_2=\epi(f)-\hypo(g)-\{(0,\alpha)\}. \end{equation*} It follows from Lemma~\ref{Qri_trans} and the imposed assumptions in (b) that the set $\Omega_1-\Omega_2$ is
quasi-regular. This allows us to apply Theorem~\ref{Theoqricap}, which ensures that the sets $\Omega_1$ and $\Omega_2$ are properly separated. Thus there exists a pair
$(\bar{u}^*,\bar{\beta})\in X^*\times\mathbb{R}$ satisfying the following two conditions: \begin{equation*}\inf_{(x,\lambda)\in\Omega_1}\big\{\la\bar{u}^*,x\ra+\bar{\beta}\lambda\big\}\ge
\sup_{(y,\mu)\in\Omega_2}\big\{\la\bar{u}^*,y\ra+\bar{\beta}\mu\big\}, \end{equation*} \begin{equation*} \sup_{(x,\lambda)\in\Omega_1}\big\{\la\bar{u}^*,x\ra
+\bar{\beta}\lambda\big\}>\inf_{(y,\mu)\in\Omega_2}\big\{\la\bar{u}^*,y\ra+\bar{\beta}\mu\big\}. \end{equation*} This gives us a constant $\gg\in\mathbb R$ such that
\begin{equation}\label{C}\la\bar{u}^*,x\ra+\bar{\beta}\lambda\ge\gg\ge\la\bar{u}^*,y\ra+\bar{\beta}\mu\end{equation} whenever $(x,\lambda)\in\Omega_1$ and $(y,\mu)\in\Omega_2$. If
$\bar{\beta}=0,$ then we have \begin{equation*}\inf_{x\in\dom(f)}\big\{\la\bar{u}^*,x\ra\big\}\ge\sup_{y\in\dom(g)}\big\{\la\bar{u}^*,y\ra\big\},\end{equation*}
\begin{equation*}\sup_{x\in\dom(f)}\big\{\la\bar{u}^*,x\ra\big\}>\inf_{y\in\dom(g)}\big\{\la\bar{u}^*,y\ra\big\}. \end{equation*} Thus the sets $\dom(f)$ and $\dom(g)$ are properly separated,
which implies by the characterization of Theorem~\ref{Theoqricap} that \begin{equation*}\qri\big(\dom(f)\big)\cap\qri\big(\dom(g)\big)=\emp\end{equation*} under the assumptions made. The obtained
contradiction verifies that $\bar{\beta}\ne 0$.

It follows from the structure of $\Omega_1$ that for any fixed $x_0\in\dom(f)$ we have $(x_0,f(x_0)+k)\in\Omega_1$ whenever $k\in\mathbb{N}$. Thus we deduce from \eqref{C} that \begin{equation*}
\la\bar{u}^*,x_0\ra+\bar{\beta}\big(f(x_0)+k\big)\ge\gg\;\mbox{ for all }\;k\in\mathbb{N}, \end{equation*} which clearly yields $\bar{\beta}>0$. It also follows from \eqref{C} that
\begin{equation*}\Big\la\frac{\bar{u}^*}{\bar{\beta}},x\Big\ra+f(x)\ge\gg\ge\Big\la\frac{\bar{u}^*}{\bar{\beta}},y\Big\ra+g(y)+\alpha\;\mbox{ for all }\;x\in\dom(f),\;y\in\dom(g). \end{equation*}
Letting $\bar{x}^*:=-\bar{u}^*/\bar{\beta}$ and $\bar{\gg}:=-\gg$ brings us to the inequalities \begin{equation}\label{eq-Inequal}\begin{array}{ll} f(x)\ge\la\bar{x}^*,x\ra-\bar{\gg}\;\mbox{ and
}\;\la\bar{x}^*,y\ra-\bar{\gg}\ge g(y)+\alpha\\\mbox{for all }\;x\in\dom(f)\;\mbox{ and all }\;y\in\dom(g).\end{array}\end{equation} The first one in \eqref{eq-Inequal} shows that
\begin{equation*}\bar{\gg}\ge\sup\big\{\la\ox^*,x\ra-f(x)\ \big|\ x\in\dom(f)\big\}=f^*(\ox^*), \end{equation*} while the second inequality in \eqref{eq-Inequal} tells us that \begin{equation*}
\bar{\gg}+\alpha\le\inf\big\{\la\ox^*,y\ra-g(y)\ \big|\ y\in\dom(g)\big\}=g_*(\ox^*).\end{equation*} Thus $\alpha\le g_*(\ox^*)-f^*(\ox^*)$, which completes the proof of the theorem.
$\h$

In the rest of this subsection we present  two useful consequences of Theorem~\ref{TheoFenchel}. The first result is based on the efficient condition for quasi-regularity via the {\em  SNC
property} from Definition~\ref{snc-conv}.

\begin{Corollary} Let $X$ be a Hilbert space, and let $f\colon X\to(-\infty,\infty]$ and $g\colon X\to[-\infty,\infty)$ be as in Theorem~{\rm\ref{TheoFenchel}}. Suppose that all the sets $\dom
(f)-\dom(g)$, $\epi(f)$, and $\epi (f)-\hypo(g)$ are closed and SNC with nonempty intrinsic relative interiors, and that the qualification condition \begin{equation*} \qri\big(\dom(f)\big)\cap
\qri\big(\dom(g)\big)\ne\emp\end{equation*} is satisfied. Then we have the Fenchel duality \eqref{duality_Theo1}. \end{Corollary}
{\bf Proof.} Theorem~\ref{quasi-reg} tells us that the sets
$\dom(f)-\dom(g)$, $\epi(f)$, and $\epi(f)-\hypo(g)$ are quasi-regular under the imposed assumptions. Applying Theorem~\ref{TheoFenchel}, we arrive at the conclusion of the corollary. $\h$

The next consequence of Theorem~\ref{TheoFenchel} involves the {\em relative interior} notion for convex sets in LCTV spaces defined in \eqref{ri}. Recall that, in contrast to the case of
finite-dimensional spaces, nonempty convex sets may have empty relative interiors in infinite dimensions.

\begin{Corollary}\label{CoroFenchelRi}  Let $X$ be an LCTV space, and $f\colon X\to(-\infty,\infty]$ and $g\colon X\to[-\infty,\infty)$ be as in Theorem~{\rm\ref{TheoFenchel}}. Suppose that the sets
$\dom(f)-\dom(g),\epi(f)$, and $\epi(f)-\hypo(g)$ have nonempty relative interiors and that the qualification condition \begin{equation}\label{ri-qc1}\ri\big(\dom(f)\big)\cap
\ri\big(\dom(g)\big)\ne\emp\end{equation} is satisfied. Then we have the Fenchel duality \eqref{duality_Theo1}.\end{Corollary}
 {\bf Proof.} Since the sets $\dom(f)-\dom(g),\epi(f)$, and
$\epi(f)-\hypo(g)$ have nonempty relative interiors, we apply
Theorem~\ref{ri-rel} and conclude that they are quasi-regular. The
duality result now follows from Theorem~\ref{TheoFenchel}.
$\h$\\[1ex]
 {\bf Acknowledgements}. The authors are very grateful to
the anonymous referee for his/her valuable remarks and suggestions
that allowed us to improve the original presentation.

\end{document}